# Further Results on Signed Product Cordial Labeling

S. Soundar Rajan[a] and J. Baskar Babujee[b]

**Abstract**
In this paper, we look into Signed Product Cordial Labeling for Splitting Graphs of Bull graph and Star graph $Spltg(K_{1,n})$, Square of Path graph $P_n^2$, Corona $C_n \odot 3k_1$ and also for the graph obtained by joining two copies of Helm $H_4$ by a Path of arbitrary length.

## 1. Introduction

A graph labeling is the assignment of integers to vertices, edges, or both if certain conditions have been met. If the domain of the mapping is the set of vertices(or edges),then the labeling is called vertex labeling(or an edge labeling). Beineke and Hegde (2001) elaborates graph labeling as a fusion of graph theory and number theory. Cahit (1987) developed the idea of cordial labeling from both graceful and harmonious labelings. Graph labeling scheme is already been applied in so many needful area like graph factorization problems. In practical situations labeling is mainly applied radar pulse code design, Neural networks, addressing systems in communication networks and in frequency assignment problem (Hale, 1980). Using graph labeling some interesting games and puzzles have been solved ,a similar work is done by Tuza (2017).

Signed product cordial labeling is a modified version of cordial labeling which was introduced by Babujee and Loganathan (2011) and they proved for Path graph, trees and cycles admits signed product cordial. In this study we investigated the existence of signed product cordial labeling of splitting graphs and some graphs with special cases.

A vertex labeling of graph G defined by $\alpha: V(G) \to \{1,-1\}$ with induced edge labeling $\alpha^*: E(G) \to \{1,-1\}$ defined by $\alpha^*(uv) = f(u)f(v)$ is called signed product cordial labeling if $|v_\alpha(-1) - v_\alpha(1)| \leq 1$ and $|e_{\alpha^*}(-1) - e_{\alpha^*}(1)| \leq 1$, where $v_\alpha(-1)$ is the number of vertices labeled with is $-1$ and $v_\alpha(1)$ is the number of vertices labeled with 1 and $e_{\alpha^*}(1)$ is the number of edges labeled with 1 and $e_{\alpha^*}(-1)$ is the number of edges labeled with $-1$. A graph G is a signed product cordial graph if it allows signed product cordial labelling.

Few definitions of graph structures are as follows. The splitting graph $Spltg(G)$ of a graph G is obtained by adding a new vertex $v'$ corresponding to each vertex $v$ of G such that $Nbhd(v) = Nbhd(v')$. An undirected planar triangle with 5 vertices is called a Bull graph. The *corona* of two graphs two graphs $G_1$(with $v_\alpha$ vertices and $e_{\alpha^*}$ edges) and $G_2$ (with $v_\beta$ vertices and $e_{\beta^*}$ edges) is defined as a graph obtained by taking one copy of $G_1$ and $n_1$ copies of $G_2$, and then joining the $i^{th}$ vertex of $G_1$ with an edge to every vertex in $i^{th}$ copy of $G_2$. Let $P_n^2$ the Square of Path graph is obtained from a path $P_n$ by joining any two vertices $u$ and $v$ by an edge if the distance between them is 2. A *Helm* $H_n$, $n \geq 3$ is the graph obtained from the wheel $W_n$ by adding a pendant edge at each vertex on the rim of the wheel $W_n$.

## 2. Main Results

**Theorem 2.1:** *The splitting graph $Spltg(K_{1,n})$ of star graph $K_{1,n}$ admits signed product cordial labeling.*

**Proof:** Let $K_{1,n}$ be a star graph with vertex set $V(G) = \{v_0, v_1, v_2, \ldots, v_n\}$ where $v_0$ is the apex vertex and $\{v_1, v_2, \ldots, v_n\}$ are pendent vertices of $K_{1,n}$ and the Edge set is defined as $E(G) = \{v_0 v_i : 1 \leq i \leq n\}$.

Now the splitting graph $Spltg(K_{1,n})$ of $K_{1,n}$ is obtained by duplicating new vertices $\{v_0', v_1', v_2', \ldots, v_n'\}$ corresponding to the vertices $\{v_0, v_1, v_2, \ldots, v_n\}$ and connect the new edges such that $Nbhd(v) = Nbhd(v')$.

The vertex set and Edge set of $Spltg(K_{1,n})$ is given by
$V(Spltg(K_{1,n})) = \{v_i : 0 \leq i \leq n\} \cup \{v_i' : 0 \leq i \leq n\}$
$E(Spltg(K_{1,n})) = \{v_0 v_i\} \cup \{v_0 v_i'\} \cup \{v_0' v_i'\}, for\ 0 \leq i \leq n$
If G is a $(V, E)$ graph then, $Spltg(G)$ is a $(V', E')$ graph such that $|V'| = 2|V|$ and $|E'| = 3|E|$.
Define the vertex labeling as
for $1 \leq i \leq n$
$\alpha(v_i) = \begin{cases} 1; i \equiv 1 (mod\ 2) \\ -1; i \equiv 0 (mod\ 2) \end{cases}$
$\alpha(v_i') = -\alpha(v_i)$
$\alpha(v_0) = 1$ and $\alpha(v_0') = -\alpha(v_0)$
The induced edge labeling is given by

[a,b]Department of Mathematics, Anna University, MIT Campus, Chennai-44, India. Email: soundarrajan@mitindia.edu, baskarjee@annauniv.edu

$\alpha^*(v_0 v_i)$
$= \begin{cases} 1; \alpha(v_0) \text{ and } \alpha(v_i) \text{ have same sign} \\ -1; \alpha(v_0) \text{ and } \alpha(v_i) \text{ have different sign} \end{cases}$
$1 \leq i \leq n$

$\alpha^*(v_0 v_i')$
$= \begin{cases} 1; \alpha(v_0) \text{ and } \alpha(v_i') \text{ have same sign} \\ -1; \alpha(v_0) \text{ and } \alpha(v_i') \text{ have different sign} \end{cases}$
$1 \leq i \leq n$

$\alpha^*(v_0' v_i')$
$= \begin{cases} 1; \alpha(v_0') \text{ and } \alpha(v_i') \text{ have same sign} \\ -1; \alpha(v_0') \text{ and } \alpha(v_i') \text{ have different sign} \end{cases}$
$1 \leq i \leq n$

The vertex and edge conditions of $Spltg(K_{1,n})$ using the above labeling pattern is given in the following table.

TABLE 1. **Vertex and edge Conditions of $Spltg(K_{1,n})$**

| $n$ | $\mathbf{v_\alpha(1)}$ | $\mathbf{v_\alpha(-1)}$ | $\mathbf{|v_\alpha(-1) - v_\alpha(1)|}$ |
|---|---|---|---|
| $n \equiv 0 \pmod 2$ | $(n+1)$ | $(n+1)$ | 0 |
| $n \equiv 1 \pmod 2$ | $(n+1)$ | $(n+1)$ | 0 |

| $n$ | $\mathbf{e_{\alpha^*}(1)}$ | $\mathbf{e_{\alpha^*}(-1)}$ | $\mathbf{|e_{\alpha^*}(-1) - e_{\alpha^*}(1)|}$ |
|---|---|---|---|
| $n \equiv 0 \pmod 2$ | $\frac{3n}{2}$ | $\frac{3n}{2}$ | 0 |
| $n \equiv 1 \pmod 2$ | $\frac{3n+1}{2}$ | $\frac{3n-1}{2}$ | 1 |

Hence $Spltg(K_{1,n})$ admits signed product cordial labeling. Signed product cordial labeling of $Spltg(K_{1,8})$ is given in Figure 1.

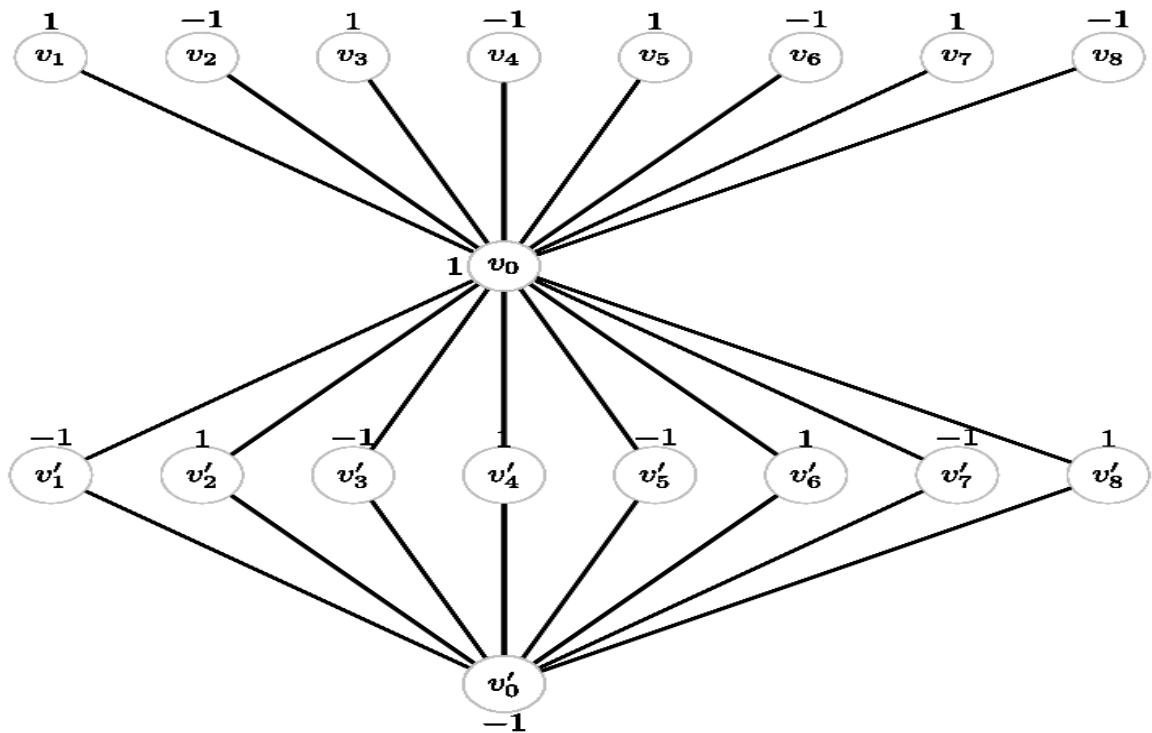

*Figure 1.* **Signed Product Cordial Labeling of $Spltg(K_{1,8})$**

***Theorem 2.2:*** *The splitting graph of Bull graph $Spltg(B_G)$ admits signed product cordial labeling.*
***Proof:***

Let $\{v_1, v_2, v_3, v_4, v_5\}$ be the vertex set of bull graph and the edge set $E(G) = \{V_i V_{i+1}\} \cup \{V_2 V_4\}$ for $1 \leq i \leq 4$. Now the splitting graph of bull graph is obtained by duplicating the new vertices $\{v_1', v_2', v_3', v_4', v_5'\}$ corresponding to the vertices $\{v_1, v_2, v_3, v_4, v_5\}$ and connect the edges satisfying $Nbhd(v) = Nbhd(v')$.

For $0 \leq i \leq 5$, Define the vertex labeling as $\alpha(v_1) = -1$

$\alpha(v_i) = \begin{cases} 1 & i \equiv 0 \pmod 2 \\ -1 & i \equiv 0 \pmod 3 \\ 1 & i \equiv 2 \pmod 3 \end{cases}$ and

$\alpha(v_i') = -\alpha(v_i)$

Now the number of vertices labeled with 1 are $v_\alpha(1) = 5$ and the number of vertices labeled with $-1$ are $v_\alpha(-1) = 5$. Therefore the total difference between the vertices labeled with 1 and $-1$ are given by $|v_\alpha(1) - v_\alpha(-1)| = 0 \leq 1$. The total number of edges labeled with 1 are $e_{\alpha^*}(1) = 8$ and the total number of edges labeled with $-1$ are $e_{\alpha^*}(-1) = 7$. Therefore the total difference between the edge labels with 1 and $-1$ are $|e_{\alpha^*}(1) - e_{\alpha^*}(-1)| = 1$. Hence the Bull graph admits signed product cordial labeling. Signed product cordial labeling of bull graph is given in Figure 2.

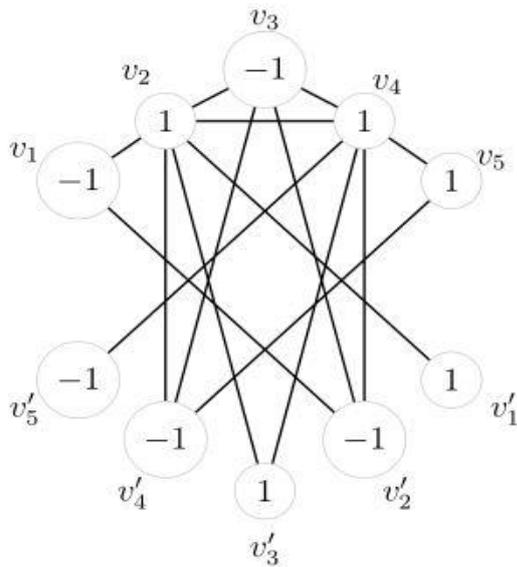

Figure 2. **Signed Product Cordial Labeling of** $Spltg(B_G)$.

**Theorem 2.3:** *The Square of path graph $P_n^2$ with $n \geq 3$ admits signed product cordial labeling.*
**Proof:** Let $P_n^2 = G(V,E)$ with $n \geq 3$ be the square of path graph with $n$ vertices and $2n-3$ edges.

Let $\{v_1, v_2, v_3 \ldots v_n\}$ be the vertices of $P_n^2$ and $\{e_1, e_2, e_3, \ldots e_n\}$ be the edges of $P_n^2$.

Define the vertex labeling for the path $P_n$ as $f: V(G) \to \{1,-1\}$ such that
$$\alpha(v_i) = \begin{cases} 1; i \text{ is odd} \\ -1; i \text{ is even} \end{cases}$$
and the induced edge labeling as

$$\alpha^*(v_i v_{i+1}) = \begin{cases} 1 \text{ both have the same sign} \\ -1; \text{both have different sign} \end{cases}$$

The induced edge labeling for $P_n^2$ is $\alpha^*(v_i v_{i+2}) = 1 \text{ for } (1 \leq i \leq n-2)$.

**Case I:** When $n$ is even

Let $v_\alpha(1) = \frac{n}{2}$ be the total number of vertices labeled with 1 and $v_\alpha(-1) = \frac{n}{2}$ be the total number of vertices labeled with $-1$. The total difference between the vertices labeled with 1 and $-1$ is $|v_\alpha(-1) - v_\alpha(1)| = 0 \leq 1$. For the edge set, $e_{\alpha^*}(1) = n-2$ be the total number of edges labeled with 1 and $e_{\alpha^*}(-1) = n-1$ be the total number of edges labeled with $-1$. The total difference between the edge labels with 1 and $-1$ are $|e_{\alpha^*}(-1) - e_{\alpha^*}(1)| = |(n-1) - (n-2)| = 1$.

**Case II:** When $n$ is odd

Let $v_\alpha(1) = \frac{n+1}{2}$ be the total number of vertices labeled with 1 and $v_\alpha(-1) = \frac{n-1}{2}$ be the total number of vertices labeled with $-1$. The total difference between the vertices labeled with 1 and $-1$ is $|v_\alpha(-1) - v_\alpha(1)| \leq 1$. Now for the edge set $e_{\alpha^*}(1) = n-2$ be the total number of edges labeled with 1 and $e_{\alpha^*}(-1) = n-1$ be the total number of edges labeled with $-1$. The total difference between the edge labels with 1 and $-1$ are $|e_{\alpha^*}(-1) - e_{\alpha^*}(1)| = |(n-1) - (n-2)| = 1$. Hence the graph $P_n^2$ with $n \geq 3$ admits signed product cordial labeling.

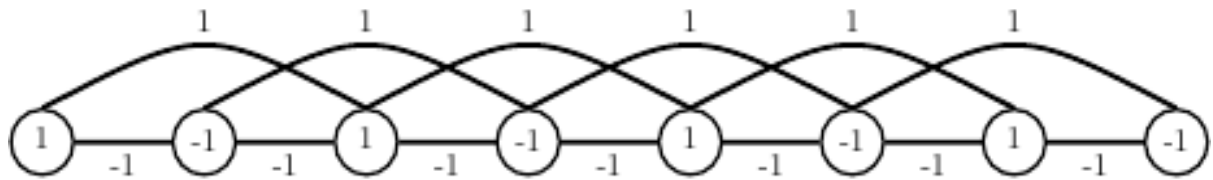

Figure 3. **Signed Product Cordial Labeling of $P_8^2$**

**Theorem 2.4:** *The corona $C_n \odot 3k_1$ admits signed product cordial labeling.*
**Proof:** Let the vertex set of $G$ be
$V(G) = \{u_x, v_x, w_x \text{ and } t_x \ 1 \leq x \leq n\}$ and the edge set of $G$ is
$E(G) = \{u_x u_{x+1}, 1 \leq x \leq n-1\} \cup \{u, u_n\} \cup \{u_x v_x, 1 \leq x \leq n\} \cup \{u_x w_x, 1 \leq x \leq n\} \cup \{u_x t_x, 1 \leq x \leq n\}$

Define the vertex labeling $\alpha: V(G) \to \{1,-1\}$ as
$u_x = 1 \quad 1 \leq x \leq n$
$v_x = -1 \quad 1 \leq x \leq n$
$w_x = 1 \quad 1 \leq x \leq n$
$t_x = -1 \quad 1 \leq x \leq n$

The induced edge labeling $\alpha^*: E(G) \to \{1,-1\}$ becomes
$\alpha^*(u_x u_{x+1}) = 1; \ 1 \leq x \leq n-1$
$\alpha^*(u_x v_x) = -1; \ 1 \leq x \leq n$
$\alpha^*(u_x w_x) = 1; \ 1 \leq x \leq n$
$\alpha^*(u_x t_x) = -1; \ 1 \leq x \leq n$
and $\alpha^*(uu_n) = 1$.

Let $v_\alpha(1) = \frac{n}{2}$ be the total number of vertices labeled with 1 and $v_\alpha(-1) = \frac{n}{2}$ be the total number of vertices labeled with $-1$. The total difference between the vertices labeled with 1 and $-1$ is $|v_\alpha(-1) - v_\alpha(1)| = 0 \leq 1$. The total number of edges labeled with $-1$ is given by $e_{\alpha^*}(-1) = \frac{n}{2}$ and the total number of edges

labeled with 1 is given by $e_{\alpha^*}(1) = \frac{n}{2}$. Therefore, the total difference between the edges labeled with 1 and $-1$ are $|e_{\alpha^*}(-1) - e_{\alpha^*}(1)| \leq 1$.

Hence the corona $C_n \odot 3k_1$ admits signed product cordial labeling.

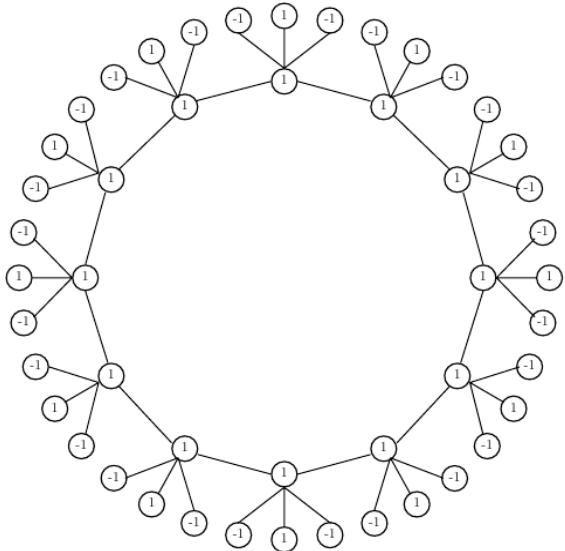

*Figure 4.* **Signed Product Cordial Labeling of $C_n \odot 3k_1$**

**Theorem 2.4:** *The two copies of $H_4$ connected to a path of arbitrary length is signed product cordial.*

**Proof :** Let $G$ be a graph obtained by joining a path $P_k$ of $(k-1)$ length and two copies of helm $H_4$. The consecutive vertices of helm $H_4$ are given by $\{v_0, v_1, \ldots, v_n, v_1', v_2', \ldots v_n'\}$, here $v_0$ is the apex vertex, $\{v_0, v_1, \ldots, v_n\}$ be the internal vertices and $\{v_1', v_2', \ldots v_n'\}$ be the external pendant vertices. Similarly let $\{w_o, w_1, \ldots, w_n, w_1', w_2', \ldots w_n'\}$ be the consecutive vertices of second copy of helm $H_4$ where $w_0$ is the apex vertex $\{w_o, w_1, \ldots, w_n\}$ be the internal vertices and $\{w_1', w_2', \cdots, w_n'\}$ be the external pendant vertices. Let $\{u_1, u_2, \ldots, u_n\}$ be the vertices of path $P_k$ with $u_1 = v_0$ and $u_k = w_0$. First label the internal vertices of first copy of $H_4$ by 1 and internal vertices of second copy of $H_4$ by $-1$. Then the external vertices of first copy of $H_4$ is labeled by $-1$ and external pendant vertices of second copy of $H_4$ is labeled by 1. Now the remaining task is to label the vertices of path $P_k$ for which we define labeling function $\alpha: V(G) \rightarrow \{1, -1\}$ as lable $u_1$ and $u_n$ as 1 and the remaining internal vertices of the path are as follow

$\alpha(u_i) = 1 \ (i \ is \ even)$ and
$\alpha(u_i) = -1 \ (i \ is \ odd)$.

Clearly the vertex and edge conditions are satisfied. Hence two copies of $H_4$ connected to a arbitrary path admits signed product cordial labeling. Signed product cordial Labeling of a graph obtained by joining two copies helm $H_4$ by path a $P_5$ is given in Figure 5.

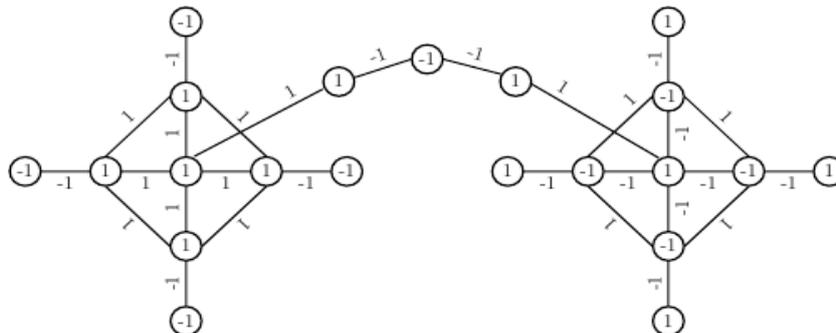

*Figure 5.* **Signed Product Cordial Labeling of Two Copies of $H_4$ by a Path $P_5$.**